\documentclass[11pt,reqno,tbtags,a4paper,final]{amsart}
\usepackage[a4paper,margin=1in,marginparwidth=2cm]{geometry}

\usepackage{amsmath,amssymb,amsfonts,amsthm,bbm}
\usepackage{xspace,xpunctuate}
\usepackage[notcite,notref]{showkeys}

\usepackage{url}
\usepackage{xcolor}
\usepackage{enumerate}
\usepackage{tikz,float}
\usepackage{subcaption,graphicx}
\usepackage{mathtools}
\usepackage[pdfencoding=auto,psdextra,breaklinks]{hyperref}
\usepackage[capitalise]{cleveref}
\usepackage{bookmark}
\usepackage{crossreftools}
\crefname{equation}{}{}
\Crefname{claim}{Claim}{Claims}

\newcommand{\ignore}[1]{\relax}

\definecolor{brown}{cmyk}{0, 0.72, 1, 0.45}
\definecolor{grey}{gray}{0.5}

\definecolor{lightRed}{cmyk}{0, 0.3, 0.3, 0.0}

\numberwithin{equation}{section}
\allowdisplaybreaks%

\theoremstyle{plain}

\theoremstyle{definition}

\theoremstyle{remark}

\DeclareMathOperator{\Prob}{\mathbb{P}} 
\renewcommand{\Pr}{\Prob}

\renewcommand{\sl}{\textit{Snakes and Ladders}\xspace}

\title[Snakes and Ladders]{%
Snakes and Ladders and Intransitivity, or
what mathematicians do in their time off}
\date{13 May 2021, rev.\ 13 Jan 2022}

\author[Gregory B. Sorkin]{Gregory B. Sorkin}
\address[Gregory B. Sorkin]{Department of Mathematics,
The London School of Economics and Political Science,
Houghton Street, London WC2A 2AE, England}
\email{g.b.sorkin@lse.ac.uk}



\begin{document}
\bibliographystyle{alpha}

\begin{abstract}
This recreational mathematics article shows that the game of \sl is intransitive:
square 69 has a winning edge over 79, which in turn beats 73, which beats 69.
Analysis of the game is a nice illustration of Markov chains,
simulations of different sorts, and size-biased sampling.
Connecting this to ``intransitive dice'' illustrates
the power of a name, and the joy of working with colleagues.
When draws do not count, we show a minimal example of intransitive dice,
with one die having just a single ``face'' and two dice each having two faces.
\end{abstract}

\maketitle

\section{Introduction}
\begin{figure}[tb]
  \centering
  \includegraphics[width=9cm]{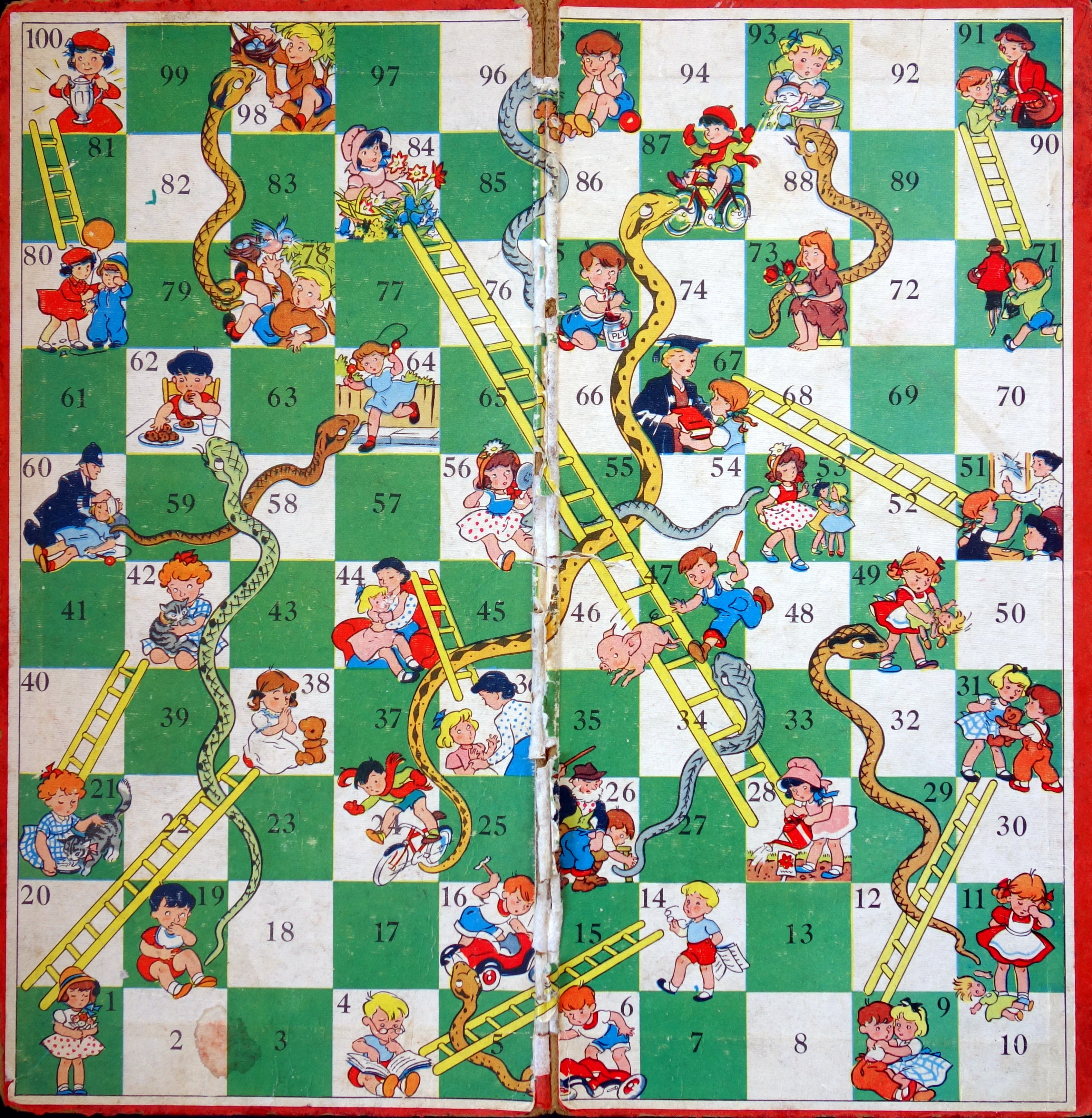}
  \caption{A sample snakes and ladders board.
  This sample, owned by the author, carried no identifying marks nor copyright.
  Squares are arranged in a serpentine pattern,
  moving to the right through the bottom row,
  left in the row above it, and so on.
  }\label{board}
\end{figure}

For those fortunate enough to be unfamiliar with it,
\sl is a children's board game of no skill and no mathematical interest.
There is a board that is in essence a strip of squares, from 1 to 100.
A sample board is shown in Figure~\ref{board}.
Each player starts in square 0 (just off the board), and players take turns.
On each turn, a player rolls a usual 6-sided die,
and advances by a number of squares equal to the roll of the die.
The twist is that some squares (``ladders''), when you land on them,
advance you to a later square,
and some (``snakes'') put you back to an earlier square.
The goal is to be first to reach square 100.

A key point is that the players are independent.
They could as well play separately, each counting how many moves they took,
and compare notes at the end:
the one with the smaller number of moves wins.
This, and the fact that no skill is involved (there are no choices to make),
are why I disparaged the game as being of no (mathematical) interest,
but actually there are interesting aspects.

In particular, as the game is played, both players tend to advance,
but there are frequent setbacks.
When you are at square $i$ and your opponent is at $j$,
it is natural to wonder who has the advantage.
We'll say that square $i$ is ``better than'' $j$
(and write ``$i>j$'')
if $i$ is more likely to win than $j$ is:
if the two players bet even odds on the outcome,
in the long run $i$ would win.
Are later squares always better?
Probably not, as it's probably better to have the possibility of a long ladder
just ahead of you than to be just past it.
Does it even make sense to ask what square is better,
or does it depend on your opponent's square?
Specifically, might the game be ``intransitive'':
is it possible that square $i$ is better than $j$,
and $j$ better than $k$,
but $k$ better than $i$,
so that $i>j>k>i$?

We will answer these questions.
We are not aware of the intransitivity question
having been asked before for \sl.
Along the way, we'll visit Markov chains,
simulation, a paradox of size-biased sampling of geometric random variables,
and intransitive dice.

First, a soup\c{c}on of history and some pesky details.

\section{History, details, and the Markov chain}

Snakes and Ladders is widely agreed to derive from an Indian game,
called \emph{gy\=an chaupar} in Hindi,
making its way to Victorian England
as a side effect of British colonialism.
An item on the \sl Wikipedia page \cite{Wsnakes}
asserts without attribution that the game has been played since the 2nd century AD,
while a number of sources credit its origin to the 13th century Sant Dnyaneshwar,
but again without citing any basis.
A scholarly article \cite{Topsfield} cites concrete evidence for the game's play
in the 18th century and says it ``is doubtless much older'',
but that since the board materials are ephemeral,
``[u]ntil earlier evidence is available,
the origins \ldots of the game must remain obscure.''
The boards vary in size, the number of snakes and ladders,
and their positions, depiction, and labelling,
but the game play remains the same.
On both continents, the game was meant to be morally educational.
Virtuous ladders, and vices represented by snakes,
would bring you towards or away from some version of heaven.
Their depiction and labelling would suit the morality of the time and place,
a Victorian version, for example,
having a ladder of Penitence leading to a square of Grace.

Whether or not morally instructive,
the game is a fine illustration of randomness.

\smallskip

An important detail for us is what it means to win.
One definition is that if you are the first to finish you win,
but that would give an unfair advantage to the first player.
(A Markov chain analysis of this version of \sl
is given in \cite{AMQ}.)
In our house we play fair:
the game goes in rounds (in each round, player 1, then player 2),
and a player wins if they finish in a round and the other player does not.
So, if player 1 finishes, player 2 has one last turn:
if they also finish, the game is a draw.
Either way, we can consider the two players separately
and simply count how many rounds it takes for each to finish:
in our fair version,
the player finishing in an earlier round wins,
and if both finish in the same round it is a draw.

With this fair version, $i>j$ means that if
one player is in square $i$ and the other in $j$,
in the same round, then $i$ wins more often than $j$
(with draws not counting either way):
in the long run, $i$ has a winning edge.

Our main interest is in ``intransitivity'', i.e.,
if there is a ``triangle'' (or longer cycle) of squares
where $i>j>k>i$.
The notion of intransitivity is natural in the fair version,
where we consider squares $i$, $j$, and $k$ all in the same round.
(It is less natural in the unfair version,
we must take into account whose move it is.
There, perhaps we'd look for a 4-cycle where
square $i$ having the move has an edge over square $j$ without it,
which in turn has an edge over square $k$ with,
that over square $l$ without, and that over square $i$ with.)

Let's return to the game's setup and clarify some details.
First, an example.
In our board, there is a ladder from square 4 to 14.
This means that square 4 can never be occupied:
if for example a player is on square 3 and rolls a 1,
they move to square 14.

\begin{figure}
  \begin{tabular}{rrrrrrrrrr}
  38 & 2 & 3 &14 & 5 & 6 & 7 & 8 &31  &10
\\11 &12 &13 &14 &15 & 6 &17 &18 &19  &20
\\42 &22 &23 &24 &25 &26 &27 &84 &29  &30
\\31 &32 &33 &34 &35 &44 &37 &38 &39  &40
\\41 &42 &43 &44 &45 &46 &26 &48 &11  &50
\\67 &52 &53 &54 &55 &53 &57 &58 &59  &60
\\61 &19 &63 &60 &65 &66 &67 &68 &69  &70
\\91 &72 &73 &74 &75 &76 &77 &78 &79 &100
\\81 &82 &83 &84 &85 &86 &24 &88 &89 & 90
\\91 &92 &73 &94 &75 &96 &97 &78 &99 &100
\\99 &78 &97 &96 &75
\end{tabular}
  \caption{Our \sl board,
  shown here laid out from top to bottom and left to right.
  Think of the values shown
  as being indexed, by position, from 1 to 105.
  Your previous value, plus your die roll (from 1 to 6),
  gives the index of your new value
  (your Markov chain state, and square number on the usual \sl board).
  For example, from the starting state of 0 (not shown), a roll of 1
  brings you to index 1, and thus to state 38,
  because square 1 has a ladder to 38.
  Squares 2 and 3 are straightforward,
  square 4 has a ladder from 4 to 14, and so forth.
  The value at index 95 is 75 because there is a snake from square 95 to square 75.
  The finishing square of 100 is followed by the squares to which
  you are redirected (your next state) if you overshoot 100.
  }\label{MCboard}
\end{figure}

There are two minor details.
One is how you finish.
If you overshoot 100, does that count as a finish,
do you stay in the same square to try again on the next turn,
or do you ``reflect'' back from 100?
We arbitrarily choose the ``reflecting'' version:
for example from 99, a roll of 3 would bring you 1 step forward to 100,
then 2 steps back to 98 \ldots where on our board there is a snake,
so you'd wind up at 78.
A second detail is that sometimes the game is played that,
if a player rolls a 6,
they are allowed an extra roll in the same turn;
it makes no essential difference, and we eschew this complication.

\medskip

To recapitulate, in essence, the game consists of a set of squares or ``states''.
From each state, there are 6 possible next states,
the actual one depending on the roll of the die.
This defines a Markov chain.
(See Figure~\ref{MCboard}.)
Our board has 84 states including 0 and 100:
squares that are the starting point of a snake or ladder do not appear as states
since it is impossible to wind up in such a square.
The winner is the first player to reach a specified state (100, in our case).

\section{Expected time to finish}

Let's return now to our questions.
It's natural to wonder, first, to what degree being further along the board
is actually helpful, and by how much.
For each state (each board square that is not the start of a snake or ladder),
what is the \emph{expected} number of moves until finishing:
the number of moves it would take, on average, over an infinite number of games?

\begin{figure}
  \centering
  \includegraphics[width=9cm]{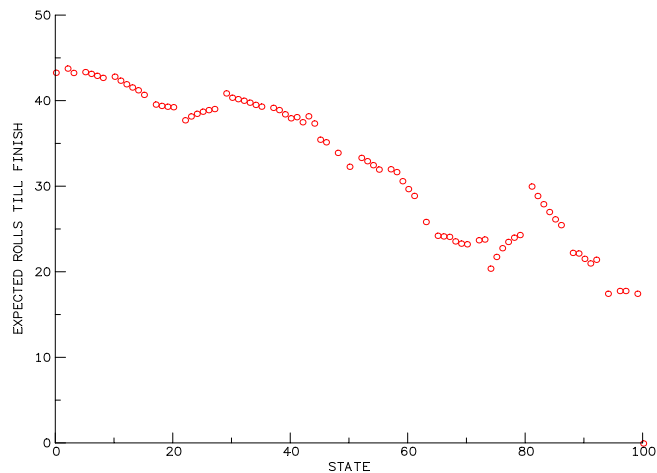}
  \caption{Expected number of moves till finish, versus state (square number).
  The game itself starts in state 0.
  The ending state is 100, so from it, the expected time to finish is 0.
  There are gaps at, for example, positions 1 and 4, because these
  boards squares are not states: in this case, they are ladders to 38 and 14.
  Being further along the board is generally helpful,
  but not consistently so.
  }\label{expectationVstate}
\end{figure}

\cref{expectationVstate} shows that indeed it is generally better to be further
advanced along the board --- later squares have a lower expected time to completion
--- but there are many exceptions.
For instance, the situation is successively \emph{worse} from squares 22 to 27,
because square 28 is a ladder to 84,
and being a bit earlier maximises the chance of landing at that ladder.

The results shown here are drawn from more detailed results giving,
for each state $i$,
the probability that, starting in state $i$,
the game finishes within $k$ moves;
in principle this should be done for all $k$ from 0 through $\infty$,
but in practice, the chance the game has not ended after 1,000 moves
is less than $10^{-14}$ (even for the worst-case starting square)
so we limited calculation to this.
This is equivalent to knowing, for each $k$, the probability that the game
ends precisely on the $k$th roll:
the differences in successive ``by time $k$'' probabilities
are the ``at time $k$'' ones,
and the cumulative sums of the ``at $k$'' probabilities are the ``by $k$'' ones.
There are two methods of going about finding this information:
simulation of the game, or calculation from the Markov transition matrix.

\subsection{Simulation}
Since these questions were just a flicker of curiosity, not a serious research agenda,
it was natural to address them by a quick and easy simulation.
We can program a computer to start a player in a specified square $i$,
perform a simulated die roll, advance the player accordingly,
and stop when the player finishes.
Repeating this for square $i$ gives a sample of the game lengths
that, with a large number of repetitions, should be an accurate
sample of the true distribution of the duration.

Observing that the game is memoryless, the simulation can be done much more efficiently.
Memorylessness means that, if we are in square $i$,
the remaining time until the end of the game
is independent of what came earlier
(though of course random depending on the future die rolls).
Thus, instead of getting just a single duration out of one game simulation,
we can get many.
Suppose a simulated game visits squares $0,2,6,10,6,9,\ldots,100$,
with 50 steps after the 0.
(From square 10, a roll of 6 brings you to square 6 via a ladder at 16.)
This play gives 51 simulated values:
from 0 the game ended in 50 steps,
from 2 in 49 steps,
and so on,
until from 100 it ended in 0 steps.
Note that from the first 6 the game ends in 48 steps,
and from the second 6 in 46 steps:
the simulation can give several remaining-time samples for a single $i$.
All in all, a play of $n$ steps gives $n$ samples (ignoring the final 100),
much better than playing a whole game to get just one sample.

\subsection{Size-biased sampling}
The latter method, if you look at it from a certain angle, appears wrong.
If we start a simulation from $i$, clearly the duration from that visit of $i$
is what we want.
The memorylessness tells us that, for later visits to $i$,
the time remaining until the finish is also a valid sample.
But, paradoxically, those later visits to $i$ obviously have
shorter game durations than the starting one:
this approach seems wrong.
In fact it is right; the mystery lies in size-biased sampling.
(The ``bus waiting time paradox'' is a beautiful example.
If buses come randomly but about once an hour,
the expected time from one bus to the next is 1 hour.
But from the moment you arrive at the bus stop,
the expected time to the next bus in an hour and, by symmetry,
the expected time since the previous bus is also an hour,
giving an expected time of 2 hours between these two buses.
This appears paradoxical.)

Here, intuitively, while it is true that looking at later visits to $i$
would lead to smaller estimates of the game duration
(certainly compared to the first visit to $i$ in the same game simulation),
countering this is that long games, with more visits to $i$,
are over-represented in the sampling.
It's not obvious that these two effects exactly balance one another,
but --- trusting to the memorylessness perspective --- they must.

To check, we can calculate.
From state $i$,
let $p$ be the probability that $i$ is visited again before the end of the game.
In that event, let $D$ be the distribution of time until the next visit to $i$.
And, let $D'$ be the distribution of time from the last visit of $i$
until the game's end.
For a given visit to $i$,
let $K$ be the number of visits to $i$ during the game
(including this visit, but no earlier ones, if this was not the first).
Conditional upon $K=k$, the length of the game is
$\sum_{s=1}^{k} A_s + B$ where $A_s \sim D$ are independent rvs
(random variables)
for the revisit durations, and $B \sim D'$ the time to get from the
final visit of $i$ to the finish.
So, $K$ tells us everything:
if the two methods of simulation result in the same distribution of $K$,
then they give (in the long run) the same sampling of game durations.
This, then, is just a question of two ways of sampling
the geometric random variable $K$.

For the first method of simulation,
$K$ is just geometrically distributed with parameter $p$:
$$ \Pr(K=k) = p^{k-1} (1-p) . $$
For the second method, of all the visits to $i$ sampled in all the games,
we wish to know what fraction of these had exactly $k$ more visits before the game end
(including this visit but no earlier ones).
For any $k \geq 1$ this is
$$ \Pr(K=k) = \frac
 {\sum_{t \geq k} p^{t-1} (1-p)}
 {\sum_{t \geq 1} t \cdot p^{t-1} (1-p)} \text{;}$$
a game with $t$ visits to $i$ occurs with probability $p^{t-1} (1-p)$,
gives one $k$th-last visit to $i$ iff $t \geq k$,
and gives $t$ visits to $i$ in all.
It is not hard to check that this expression simplifies to $p^{k-1} (1-p)$.
That is, the fraction of $i$-visits that are $k$th-last ones in the second simulation approach
is the same as the first approach's probability that there are $k$ visits to $i$,
and the two approaches do (as they must) lead to the same result.

\subsection{Markov chain}
The \sl Markov chain, like any other,
is completely described by its transition matrix $A$.
For states $i$ and $j$, $A_{ij}$ is the probability of moving from state $i$ to state $j$
in one step.
Here, for example, $A_{17,19}=\frac16$: a die roll of 2 (only) brings us from 17 to 19.
To get from $i$ to $j$ in exactly two steps means moving from $i$ to some $k$
in one step and $k$ to $j$ in the next, which happens with probability
$\sum_k A_{ik} A_{jk} = (A^2)_{ij}$.
Repeating this gives a fundamental property of Markov processes,
that the probability of getting from $i$ to $j$ in exactly $s$ steps is
$(A^s)_{ij}$.

If we are interested in the probability,
starting from $i$, of reaching the final state 100 in $s$ steps,
here that is given by $(A^s)_{i, 100}$.
Specifically, the finishing state is ``absorbing'':
from state 100 there is probability 1 of returning to 100 ($A_{100,100}=1$)
and probability 0 of moving to any other state.
In this case $(A^s)_{i, 100}$ represents the probability of being
in the finish state at time $s$ (perhaps having reached the state earlier).
As remarked earlier, $(A^s)_{i, 100} - (A^{s-1})_{i, 100}$
is the probability that the game duration, from $i$, is exactly $s$.

So, repeating for say $s$ from 0 to $1,000$ gives the probability
$(A^s)_{i, 100}$ that, for each start state $i$, the game ends at time $s$.
In practice, this gives, for each $i$, the distribution of game lengths
(the only error being the $<10^{-14}$ fraction
of games that are longer than $1,000$ rolls).

Let's quickly return to the expected game durations from each state.
Let $f_i(s)$ be the probability that the game, starting from $i$,
ends in exactly $s$ steps.
Then, the expected duration of the game, from $i$,
is simply $$\sum_{s=0}^{\infty} s \cdot f_i(s) . $$
This leads to the results shown in Figure~\ref{expectationVstate}.

\section{Pair competitions and intransitivity}

What about the probability that a player in state $i$
finishes in fewer rounds than an opponent in state $j$?
For any state $i$, define
$ g_i(s) = \sum_{t \leq s} f_i(s) $;
this is the probability that the game has finished \emph{within} time $s$,
starting from $i$.
For $i$ to beat $j$ means that $i$ finishes in some round $s$
by which $j$ has not yet finished, so $i$ beats $j$ with probability
$$ Q_{ij} := \sum_{s=0}^{\infty} f_i(s) (1-g_j(s)) . $$
Truncating this to a finite sum
gives our estimate of the probability $Q_{ij}$ that $i$ beats $j$.
We compute this for all pairs $i,j$.
Specifically, for each $s$ we compute the array of all $f_i(s) (1-g_j(s))$
(an ``outer product'' of the vector of all $f_i(s)$ with that of all $g_j(s)$).
This calculation is perhaps not as elegant as could be,
but since the outer product array
can be computed faster than a single matrix product $A \cdot A^s$
(they are of the same dimension, square over the number of states),
it is efficient enough.

Define the ``excess'' $X_{ij}$ of $i$ over $j$ by
$X_{ij}= Q_{ij} - Q_{ji}$,
the win probability of $i$ over $j$ versus that of $j$ over $i$.
If on each game the winner got \pounds 1, with no money exchanged for a draw,
$X_{ij}$ would be $i$'s average winning, playing against $j$.
We won't need it, but the probability of a draw is just $1-Q_{ij}-Q_{ji}$.

Our original question translates to whether there are states $i,j,k$
such that $X_{ij}$, $X_{jk}$, and $X_{ki}$ are all positive.
Specifically, let's look for such states where
$X_{ij} \geq c$, $X_{jk} \geq c$, and $X_{ki} \geq c$,
and $c$ is as large as possible.
While this was done using a trick something like a matrix product,
that way is no faster than trying all triples of states,
so let us not explain but just assert that we found the best triple.

The result is that states 69, 79, and 73 form such a triangle,
each with a winning edge at least $\frac12\%$ over the next in the cycle.
Specifically, state 69 has a winning edge over state 79 of $X_{69,79} \approx 0.0077$,
with larger wins of
$X_{79,73} \approx 0.0112$ and $X_{73,69} \approx 0.0171$.
The respective win probabilities are
$Q_{69,79} \approx 0.4970$,
$Q_{79,73} \approx .4990$ and $Q_{73,69} \approx 0.4930$.

\subsection{Check by simulation}
These figures were checked against simulations.
Starting from 0, $100,000$ games were simulated,
with a total of about 4.4M die rolls;
the maximum game length observed was under 500 rolls.
Each of the three states in question was visited at least 25,000 times.
Comparing with the calculated winning edges above, namely
0.0077, 0.0112, and 0.0171,
the simulated ones were about
0.0090, 0.0096, and 0.0172.
A second simulation gave similar results:
0.0081, 0.0127, 0.0146.
Each simulation's time is dominated by the gameplay,
taking under 3 minutes in an inefficient implementation.
Simulations with 10,000 games are even less accurate,
often showing negative values rather than the positive ones desired.
A simulation with 1M games takes 5 minutes (in a quickly improved implementation)
and gives results wonderfully close to the calculated ones:
0.0072, 0.0115, and 0.0171 (each off by at most 0.0005 from the calculation).

This simulation method has the advantage that it generates data for all squares.
However, it has some implementation complications and thus possible errors.

\subsection{\ldots and check it twice}
\begin{figure}
  \centering
  \includegraphics[width=9cm]{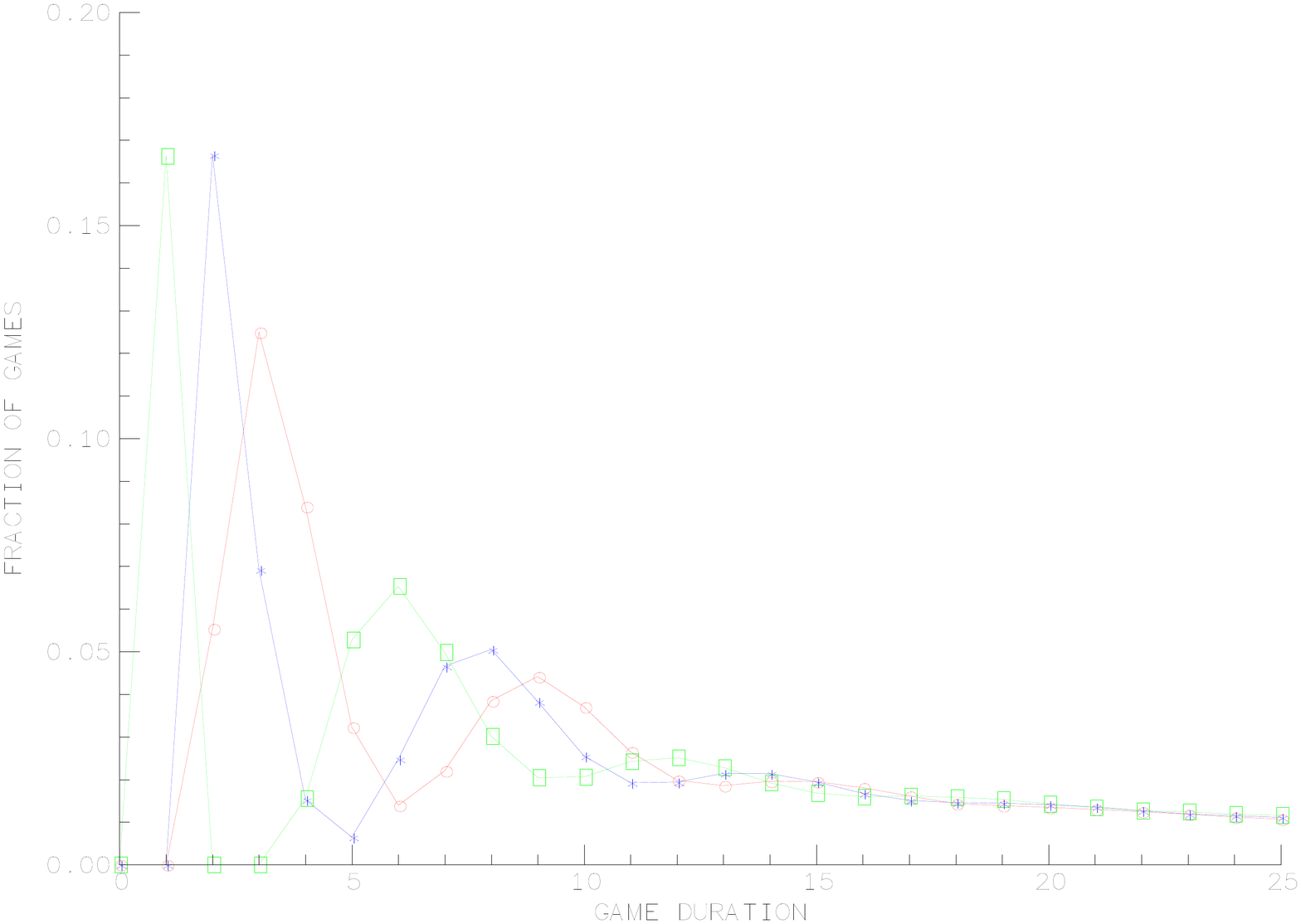}
  \caption{Fraction of simulated games of each duration,
  starting from square 69 (red), 79 (green), and 73 (blue).
  All are truncated to durations at most 25 for visual clarity;
  the maximum observed over 1M simulated games each was under 500.}\label{ijkTimes}
\end{figure}

To be sure, I also quickly tried the very simplest ``dumb'' simulation:
simulating just the time of a single game started from a given state.
I simulated 1M games from each of these three states,
and measured in what fraction of them player 1 beat player 2.
This gave an estimate of
0.0067, 0.0100, 0.0177.
The estimated time-to-finish distributions, for these squares,
from this simulation method, are shown in \cref{ijkTimes}.

Other than simplicity, an advantage of this is that we are assured that
each square is represented by 1M samples.
In the method where 1M games were played each contributing as many samples
as its length, altogether some 44M samples were generated,
but that works out to only about 500,000 per state,
and square 69 was under-represented with only about 261,000 samples.

On the other hand, the ``dumb'' method's comparison of the starting squares
``game by game'' is inefficient,
as we could just as well compare square 69's game 1 with square 79's game 17
--- in principle 1M times 1M comparisons for each pair of starting squares.
This can be done efficiently, as the game length was always under 500,
so for each starting square, a histogram of the 1M duration samples is quite compact.
This gave a better estimate,
0.0074, 0.0111, 0.0176,
using exactly the same simulation data.
The cruder method does have one advantage though:
each of the $10^6$ ``games between'' player 1 and player 2 is independent,
so it is straightforward to estimate the variance in the winning-edge estimate.
The cleverer method here makes $10^{12}$ comparisons from these $10^6$ games,
but they are not independent,
and the variance of its estimate, while smaller,
is harder to estimate.

\section{Intransitive dice}
Despite having been open to the theoretical possibility
of there being intransitivity in \sl,
I was struck by this result.
I did not know anything else like it.
I expected that something related must be known
and did an Internet search but,
not knowing the right keywords,
did not quickly find anything.
I reached out to a couple of colleagues and got an almost immediate response
from Bernhard von Stengel:
``intransitive dice''.
I will not attempt to summarise the literature,
but there is (as usual) a nice Wikipedia article \cite{Wdice}.

One example taken from the article is three six-sided dice,
die $A$ labelled $2,2,6,6,7,7$,
$B$ labelled $1,1,5,5,9,9$,
and $C$ labelled $3,3,4,4,8,8$.
Aside from the conventionality and physical practicality of six sides,
we can as well think of these as three-sided dice,
die $A$ labelled $2,6,7$,
$B$ labelled $1,5,9$,
and $C$ labelled $3,4,8$.
Note that there are no ties.
Die $A$ beats $B$ $5/9$ths of the time,
for an edge of $1/9$th,
where here winning means having a larger value.

(If we think of the die values as analogous to the durations of our game,
and desire that a smaller value should win,
we can simply replace these values with their difference from 10;
this has the effect of swapping the identities of $A$ and $C$.)

Is there an even smaller example?
This example used 9 distinct numerical values over the 3 faces of 3 dice.
Trying to use 8 values over 3 fair dice having 3, 3, and 2 faces
we find another example:
$A=4,5$, $B=1,6,7$, and $C=2,3,8$.
$A$ beats $B$ with probability (w.p.) $2/3$ (if $B$ comes up 6 or 7),
$B$ beats $C$ w.p.\ $5/9$ (if $B$ comes up 1, or $C$ is 8),
and $C$ beats $A$ w.p.\ $2/3$ (if $C$ comes up 2 or 3).

If we allow ties and only insist on a winning edge
(not necessarily winning w.p.\ greater than $1/2$)
there is an example, again with fair dice with 2, 3, and 3 faces,
but using only 7 distinct values:
$A=3,4$, $B=1,4,6$, $C=2,3,7$.
Here, $A$ beats $B$ w.p.\ $1/2$, but draws w.p.\ $1/6$th
and loses w.p.\ $1/3$, for a winning edge.
As before, $B$ beats $C$ w.p.\ $5/9$ and there are no draws.
$C$ beats $A$ symmetrically to how $A$ beats $B$.

If we allow ``unfair'' dice, with arbitrary probabilities,
a smallest example has dice with 1, 2, and 2 sides.
We can take $A=2$ (deterministically),
$B=1,3$ with probabilities $1/3, 2/3$ respectively, and
$C=1,4$ with probabilities $5/9,4/9$.
Clearly $A$ beats $B$ w.p.\ $2/3$
and $C$ beats $A$ w.p.\ $5/9$.
Finally $B$ beats $C$ w.p.\ $4/9$ (when $C=4$)
while $C$ beats $B$ w.p.\ $5/9 \times 2/3 = 10/27$ (only when $C=1$ and $B=3$),
giving this case a winning edge of $2/27$.

There can be no smaller example:
it would have to have dice with 1, 1, and 2 sides.
Without loss of generality,
one-sided die $A$ must have a smaller value than one-sided die $B$,
in which case, as $B$ wins over $C$,
$A$ must also win over $C$
(with $A$ and $B$ deterministic, the outcomes are determined by $C$ alone,
and every case for $B$ vs $C$ is at least as good for $A$ vs $C$).

\bibliography{snakes}
\end{document}